# FENCE METHODS FOR MIXED MODEL SELECTION

By Jiming Jiang,[1] J. Sunil Rao,[2] Zhonghua Gu and Thuan Nguyen

*University of California, Davis, Case Western Reserve University, ALZA Corporation and University of California, Davis*

Many model search strategies involve trading off model fit with model complexity in a penalized goodness of fit measure. Asymptotic properties for these types of procedures in settings like linear regression and ARMA time series have been studied, but these do not naturally extend to nonstandard situations such as mixed effects models, where simple definition of the sample size is not meaningful. This paper introduces a new class of strategies, known as fence methods, for mixed model selection, which includes linear and generalized linear mixed models. The idea involves a procedure to isolate a subgroup of what are known as correct models (of which the optimal model is a member). This is accomplished by constructing a statistical *fence*, or barrier, to carefully eliminate incorrect models. Once the fence is constructed, the optimal model is selected from among those within the fence according to a criterion which can be made flexible. In addition, we propose two variations of the fence. The first is a stepwise procedure to handle situations of many predictors; the second is an adaptive approach for choosing a tuning constant. We give sufficient conditions for consistency of fence and its variations, a desirable property for a good model selection procedure. The methods are illustrated through simulation studies and real data analysis.

**1. Introduction.** On the morning of March 16, 1971, Hirotugu Akaike, as he was taking a seat on a commuter train, came out with the idea of a connection between the relative Kullback–Leibler discrepancy and the empirical log-likelihood function, a procedure that was later named Akaike's

Received March 2007; revised June 2007.
[1]Supported in part by NSF Grants DMS-02-03676 and DMS-04-02824.
[2]Supported in part by NSF Grants DMS-02-03724, DMS-04-05072 and NIH Grant K25-CA89868.

information criterion, or AIC (Akaike [1, 2]; see Bozdogan [5] for the historical note). The idea has allowed major advances in model selection and related fields (e.g., de Leeuw [7]).

The procedure essentially amounts to minimizing a criterion function of the following form:

$$\hat{D}_M + \lambda_n |M|, \tag{1}$$

where $M$ represents a candidate model, $\hat{D}_M$ is a measure of lack of fit by $M$ and $|M|$ denotes the dimension of $M$, usually in terms of the number of estimated parameters under $M$. The main difference between procedures is made by $\lambda_n$, where $n$ is the sample size. This is called a "penalizer," although some authors refer $\lambda_n |M|$ as the penalizer. For example, $\lambda_n = 2$ for AIC. A number of similar criteria have since been proposed, including the Bayesian information criterion (BIC; Schwarz [26]) in which $\lambda_n = \log(n)$, and a criterion due to Hannan and Quinn (HQ; Hannan and Quinn [10]) in which $\lambda_n = c \log\{\log(n)\}$ and $c$ is a constant $> 2$. All these procedures can be viewed as special cases of the generalized information criterion (GIC; Nishii [21], Shibata [27]). A nice monograph on model selection from various perspectives is edited by Lahiri [18].

Although these criteria are widely used, difficulties are often encountered, especially in some nonconventional situations. A broad class of such nonconventional cases are mixed effects models, including linear and generalized linear mixed models. For example, consider the following linear mixed model, $y_{ij} = x'_{ij}\beta + u_i + v_j + e_{ij}$, $i = 1, \ldots, m_1$, $j = 1, \ldots, m_2$, where $x_{ij}$ is a vector of known covariates, $\beta$ is a vector of unknown regression coefficients (the fixed effects), $u_i$, $v_j$ are random effects and $e_{ij}$ is an additional error. It is assumed that $u_i$'s, $v_j$'s and $e_{ij}$'s are independent, and that for the moment, $u_i \sim N(0, \sigma_u^2)$, $v_j \sim N(0, \sigma_v^2)$, $e_{ij} \sim N(0, \sigma_e^2)$. It is well known (e.g., Hartley and Rao [11], Harville [13], Miller [20]) that in this case, the effective sample size for estimating $\sigma_u^2$ and $\sigma_v^2$ is not the total sample size $m_1 \cdot m_2$, but $m_1$ and $m_2$, respectively. Now suppose that one wishes to select the fixed covariates, which are components of $x_{ij}$ under the assumed model structure using BIC. It is not clear what should be in place of $n$ in (1), where $\lambda_n = \log(n)$ (it does not make sense to let $n = m_1 \cdot m_2$). In fact, in cases of correlated observations, such as the example here, the definition of "sample size" is often unclear.

Furthermore, suppose that normality is not assumed in the above linear mixed model. In fact, the only distributional assumptions are that the random effects and errors are independent, and that they have means zero and variances $\sigma_u^2$, $\sigma_v^2$ and $\sigma_e^2$, respectively. Now suppose that one, again, wishes to select the fixed covariates using AIC, BIC or HQ. It is not clear how to do this because all three require the likelihood function in order to evaluate $\hat{D}_M$.

In a way, model selection and estimation are two components of a process called model identification. While there is extensive literature on parameter estimation in linear and generalized linear mixed models, the other component, that is, mixed model selection, has received much less attention. Only recently have some results emerged in the area of linear mixed model selection. Datta and Lahiri [6] discussed a model selection method based on computation of the frequentist's Bayes factor in choosing between a fixed effects model and a random effects model. They focused on a one-way random effects model and noted a connection between the choice of fixed or random effects models and test of the hypothesis that the variance of the random effects is zero. Note that, however, not all model selection problems can be formulated as hypothesis testing. Jiang and Rao [16] developed various GICs suitable for linear mixed model selection and proved consistency of their procedures. Meza and Lahiri [19] demonstrated the limitations of Mallows' $C_p$ statistic in selecting the fixed covariates in a nested error regression model which is a special case of the linear mixed models. They showed by simulation results that the $C_p$ method without modification does not work well when the variance of the random effects is large; on the other hand, a modified $C_p$ criterion obtained by adjusting the intra-cluster correlations performs similarly as the $C_p$ in regression settings. Fabrizi and Lahiri [9] developed a robust model selection method in the context of complex surveys. Another related paper is Vaida and Blanchard [28], in which the authors proposed a conditional AIC where the penalty term is related to the effective degrees of freedom for a linear mixed model proposed by Hodges and Sargent [14].

It should be pointed out that all these studies are limited to linear mixed models, while model selection in generalized linear mixed models (GLMMs) has never been seriously addressed in the literature. It is well known that the likelihood function under a GLMM may involve high-dimensional integrals which are difficult to evaluate, which makes a procedure based on (1) computationally unattractive. Furthermore, our simulation results suggested that in the case of GLMM selection, a GIC procedure is much more sensitive to the choice of $\lambda_n$ than in linear mixed model selection.

In summary, the major concerns regarding the GIC procedures when applied to mixed model selection are: (i) they depend on the effective sample size which is unclear in typical situations of mixed models; (ii) they rely on the likelihood function which may not be available; (iii) they do not seem applicable to GLMMs; and (iv) their finite sample performance may be sensitive to different choices of penalties. These motivate the development of a new procedure for mixed model selection, called the *fence* method, which we describe in detail in the next section. In Section 3, we propose two variations of the fence method. The first is a stepwise fence procedure; the second is an adaptive fence procedure. In Section 4, we address the issue of consistency of different fence methods. In Section 5, we present some examples of simulations and real data analysis. Some concluding remarks are made in Section 6. Proofs of the main results are given in Section 7.

**2. The fence method.** It is illustrative to first consider a simple example. Suppose that the observations $y_{ij}$ satisfy the following linear mixed model,

$$(2) \qquad y_{ij} = x'_{ij}\beta + \alpha_i + \varepsilon_{ij},$$

$i = 1, \ldots, m$, $j = 1, \ldots, n_i$, where $x_{ij}$ is a vector of covariates, $\beta$ is a vector of unknown regression coefficients, $\alpha_i$ is a random effect and $\varepsilon_{ij}$ is an error. It is assumed that the random effects and errors are independent such that $\mathrm{E}(\alpha_i) = 0$, $\mathrm{var}(\alpha_i) = \sigma^2$, $\mathrm{E}(\varepsilon_{ij}) = 0$ and $\mathrm{var}(\varepsilon_{ij}) = \tau^2$. Even for this simple model, there are various model selection problems. For example, the selection of the fixed covariates; whether or not to include the random effects, etc. Our strategy is based on a quantity $Q_M = Q_M(y, \theta_M)$, where $y$ represents the vector of observations, $M$ indicates a candidate model and $\theta_M$ denotes the vector of parameters under $M$. It is required that $\mathrm{E}(Q_M)$ is minimized when $M$ is a true model and $\theta_M$ the true parameter vector under $M$. This means that $Q_M$ is a measure of lack-of-fit. Here by true model, we mean that $M$ is a correct model, but not necessarily the most efficient one. For example, suppose that $y_{ij}$ satisfy (2) with $x'_{ij}\beta = \beta_0 + \beta_1 x_{1ij} + \beta_2 x_{2ij}$, where all the $\beta$'s are nonzero. Then for the problem of selecting the fixed covariates, this model is optimal in the sense that the number of fixed covariates cannot be further reduced. However, the model remains true if in (2) $x'_{ij}\beta = \beta_0 + \beta_1 x_{1ij} + \beta_2 x_{2ij} + \beta_3 x_{1ij}^2$ (with $\beta_3 = 0$). But the latter model is not optimal. On the other hand, the model with $x'_{ij}\beta = \beta_0 + \beta_1 x_{1ij}$ in (2) is an incorrect model. In this paper, we use the terms "true model" and "correct model" interchangeably. Below are some options for $Q_M$ under different situations:

1. *Maximum likelihood* (*ML*) *model selection.* If the normality assumption is made regarding the random effects and errors, an example of $Q_M$ is the negative of the log-likelihood under $M$.

2. *Mean and variance/covariance* (*MVC*) *model selection.* Suppose that the situation is a bit more complicated. First, the errors are correlated within the clusters with some (parametric) covariance structure. Second, the normality assumption is not made. In such a case, the likelihood function is not available. However, one may consider $Q_M = |(T'V_M^{-1}T)^{-1}T'V_M^{-1}(y - \mu_M)|^2$, where $\mu_M$ and $V_M$ are the mean vector and covariance matrix under $M$, and $T$ is a (not necessarily square) matrix of full rank. Note that in this case, $\mu_M = X_M \beta_M$, where $X_M$ is the matrix of covariates under $M$ and $\beta_M$ the vector of regression coefficients under $M$. Thus, such a $Q_M$ may be used to select the fixed covariates as well as the (parametric) covariance structure.

A special case of the MVC is least squares (LS) model selection, in which $T = I$, the identity matrix, hence $Q_M = |y - X_M \beta_M|^2$. The latter is useful, for example, if only the fixed covariates are subject to selection while the covariance structure of the data is unknown.

It is easy to verify that all the $Q_M$ above satisfy the basic requirement of lack-of-fit above. Other choices of $Q_M$ are considered in Jiang et al. [17].

2.1. *Building the fence.* Given a specific $Q_M$, let $\hat{Q}_M = Q_M(y, \hat{\theta}_M)$, where $\hat{\theta}_M$ is the minimizer of $Q_M$ over $\theta_M \in \Theta_M$, the parameter space under $M$, that is, $\hat{Q}_M = \inf_{\theta_M \in \Theta_M} Q_M(\theta_M, y)$. Let $\tilde{M} \in \mathcal{M}$ be such that $\hat{Q}_{\tilde{M}} = \min_{M \in \mathcal{M}} \hat{Q}_M$, where $\mathcal{M}$ represents the set of candidate models. We assume that $\mathcal{M}$ contains a true model. Note that in many cases, $\tilde{M}$ can be determined without any calculation. For example, if $\mathcal{M}$ contains a full model, say $M_{\rm f}$, that is, a model such that all other models in $\mathcal{M}$ are submodels of $M_{\rm f}$, then clearly, $\tilde{M} = M_{\rm f}$ and, since $\mathcal{M}$ contains a true model, $M_{\rm f}$ is also a true model. In general, $\mathcal{M}$ may not contain a full model, but the following lemma shows that at least in large sample, $\tilde{M}$ is expected to be a correct model.

LEMMA 1. *Under Assumptions* A1–A5 *in Section 4, we have with probability tending to one that $\tilde{M}$ is a true model.*

The proof of Lemma 1 follows directly from that of Theorem 1 in the sequel.

However, the main question is, "Are there other correct models with smaller dimension than $\tilde{M}$?" This is where the fence idea comes in. As mentioned, the idea is to construct a statistical barrier, called the fence, to carefully eliminate incorrect models. Then for the models within the fence which are considered correct, one may use whatever criterion of optimality to select the optimal model. In many cases, the criterion of optimality is minimal dimension of the model, but it may be replaced by some other considerations, or incorporate scientific or economical concerns. For example, in small area estimation (SAE, e.g., Rao [24]) a main problem of interest is the prediction of small area means. Thus, some measure of the prediction errors, such as the mean squared prediction error, should be taken into account in selecting the optimal model within the fence. By the way, the linear mixed model (2), also known as the nested-error regression model (e.g., Battese, Harter and Fuller [4]), has extensive applications in SAE. The fence is constructed through the following inequality:

$$\hat{Q}_M \leq \hat{Q}_{\tilde{M}} + c_n \hat{\sigma}_{M,\tilde{M}}, \qquad (3)$$

where $\hat{\sigma}_{M,\tilde{M}}$ is an estimate of the standard deviation of $\hat{Q}_M - \hat{Q}_{\tilde{M}}$, denoted by $\sigma_{M,\tilde{M}}$. It can be shown that the latter is an appropriate measure of $\hat{Q}_M - \hat{Q}_{\tilde{M}}$ for a correct model $M$; while for an incorrect model, the difference is expected to be much larger. Furthermore, $c_n$ denotes a tuning constant. For consistency of the model selection (see Section 4), it is required that $c_n$ increase (slowly) with the sample size. Here consistency is in the sense that as the sample size increases, the probability that the procedure selects an optimal model approaches one. In Section 3.2, we show how to choose $c_n$ adaptively in order to improve the finite sample performance.

In case the minimal dimension criterion is used, an effective algorithm is outlined below.

2.2. *The fence algorithm.* For simplicity, consider the case that $\tilde{M}$ is unique. Let $d_1 < d_2 < \cdots < d_L$ be all the different dimensions of the models $M \in \mathcal{M}$. We proceed as follows:

(i) Consider $\mathcal{M}_1 = \{M \in \mathcal{M} : |M| = d_1 \text{ and } (3) \text{ holds}\}$; if $\mathcal{M}_1 \neq \varnothing$, stop (no need for any more computation). Let $M_0 \in \mathcal{M}_1$ be such that $\hat{Q}_{M_0} = \min_{M \in \mathcal{M}_1} \hat{Q}_M$; $M_0$ is the selected model.

(ii) If $\mathcal{M}_1 = \varnothing$, consider $\mathcal{M}_2 = \{M \in \mathcal{M} : |M| = d_2 \text{ and } (3) \text{ holds}\}$; if $\mathcal{M}_2 \neq \varnothing$, stop. Let $M_0 \in \mathcal{M}_2$ be such that $\hat{Q}_{M_0} = \min_{M \in \mathcal{M}_2} \hat{Q}_M$; $M_0$ is the selected model.

(iii) Continue until the program stops (it will at some point).

In short, the algorithm may be described as follows: Check the candidate models, from the simplest to the most complex. Once one has discovered a model that falls within the fence and checked all the other models of the *same simplicity* (for membership within the fence), one stops.

2.3. *Estimation of $\sigma_{M,\tilde{M}}$.* In some cases, this is straightforward. For example, suppose that the likelihood function is available, and $Q_M$ is chosen as the negative log-likelihood. Furthermore, suppose that $M_f \in \mathcal{M}$. Then under some regularity conditions $2(\hat{Q}_M - \hat{Q}_{M_f})$ has an asymptotic $\chi_d^2$ distribution with $d = |M_f| - |M|$. Thus, if $\tilde{M} = M_f$, we have $\sigma_{M,M_f} \approx \sqrt{(|M_f| - |M|)/2}$.

However, such an asymptotic $\chi^2$ distribution may not exist in general. Nevertheless, suppose that $M^*$ is true. Then in the case of clustered observations one can approximate, under some regularity conditions, $\sigma_{M,M^*}^2$ by $\mathrm{var}(Q_M - Q_{M^*})$. Furthermore, suppose that $Q_M$ can be expressed as $\sum_{i=1}^m Q_{M,i}$, where $Q_{M,i} = Q_{M,i}(y_i, \theta_M)$. Then $\mathrm{var}(Q_M - Q_{M^*}) = \mathrm{E}[\sum_{i=1}^m (Q_{M,i} - Q_{M^*,i})^2 - \sum_{i=1}^m \{\mathrm{E}(Q_{M,i}) - \mathrm{E}(Q_{M^*,i})\}^2]$. Thus, an observed variance is obtained by removing the outside expectation and replacing the parameters and inside expectations by their estimators. See Jiang et al. [17] for more detail. The latter also considered several cases of nonclustered observations, including Gaussian mixed models, non-Gaussian linear mixed models and GLMMs.

**3. Variations.** In this section, we propose two variations of the fence. The first aims at making the fence procedure computationally more attractive. The second focuses on choosing the tuning constant $c_n$ to improve the finite sample performance of the fence.

3.1. *A stepwise fence procedure.* As mentioned, the fence has the computational advantage that it starts with the simplest models, and, therefore, may not need to search the entire model space in order to determine the optimal model. On the other hand, such a procedure may still involve a lot of evaluations when the model space is large. For example, in quantitative trait loci mapping, variance components arising from the trait genes, polygenic and environmental effects are often used to model the covariance structure of the phenotypes given the identity by descent sharing matrix (e.g., Almasy and Blangero [3]). Such a model is usually complex due to the large number of putative trait loci. To make the fence procedure computationally more attractive to large and complex models, we propose the following variation of fence for situations of complex models with many predictors.

To be more specific, consider the extended GLMMs introduced by Jiang and Zhang [15]. It is assumed that given a vector $\alpha$ of random effects, the responses $y_1, \ldots, y_n$ are conditionally independent, such that $\mathrm{E}(y_i|\alpha) = h(x_i'\beta + z_i'\alpha)$, $1 \le i \le n$, where $h(\cdot)$ is a known function, $\beta$ is a vector of unknown fixed effects and $x_i$, $z_i$ are known vectors. Furthermore, it is assumed that $\alpha \sim N(0, \Sigma)$, where the covariance matrix $\Sigma$ depends on a vector $\psi$ of variance components. Let $\beta_M$ and $\psi_M$ denote $\beta$ and $\psi$ under $M$, and $g_{M,i}(\beta_M, \psi_M) = \mathrm{E}\{h_M(x_i'\beta_M + z_i'\Sigma_M^{1/2}\xi)\}$, where $h_M$ is the function $h$ under $M$, $\Sigma_M$ is the covariance matrix under $M$ evaluated at $\psi_M$, and the expectation is taken with respect to $\xi \sim N(0, I_m)$ (which does not depend on $M$). Here $m$ is the dimension of $\alpha$ and $I_m$ the $m$-dimensional identity matrix. Let $Q_M = \sum_{i=1}^{n}\{y_i - g_{M,i}(\beta_M, \psi_M)\}^2$.

Write $X = (x_i')_{1 \le i \le n}$ and $Z = (z_i')_{1 \le i \le n}$. We assume that there is a collection of covariate vectors $X_1, \ldots, X_K$, from which the columns of $X$ are to be selected. Furthermore, we assume that there is a collection of matrices $Z_1, \ldots, Z_L$ such that $Z\alpha = \sum_{s \in S} Z_s \alpha_s$, where $S \subset \{1, \ldots, L\}$, and each $\alpha_s$ is a vector of i.i.d. random effects with mean 0 and variance $\sigma_s^2$. The subset $S$ is subject to selection. The parameters under an extended GLMM are the fixed effects and variances of the random effects. Note that in this case, the full model corresponding to $X\beta + Z\alpha = \sum_{k=1}^{K} X_k\beta_k + \sum_{l=1}^{L} Z_l\alpha_l$ is among the candidate models. Thus, we let $\tilde{M} = M_{\mathrm{f}}$. The idea is to use a forward–backward procedure to generate a sequence of candidate models, among which the optimal model is selected using the fence method. We begin with a forward procedure. Let $M_1$ be the model that minimizes $\hat{Q}_M$ among all models with a single parameter; if $M_1$ is within the fence, stop the forward procedure; otherwise, let $M_2$ be the model that minimizes $\hat{Q}_M$ among all models that add one more parameter to $M_1$; if $M_2$ is within the fence, stop the forward procedure; and so on. The forward procedure stops when the first model is discovered within the fence. The procedure is then followed by a backward elimination. Let $M_k$ be the final model of the forward procedure. If no submodel of $M_k$ with one less parameter is within the fence, $M_k$ will be our selection; otherwise, $M_k$ is replaced by $M_{k+1}$ which is a submodel of $M_k$ with one less parameter and is within the fence, and so on. We call such a variation of fence the forward–backward (F-B) fence.

3.2. *Adaptive fence procedure.* In this section, we address the issue regarding choosing the tuning constant $c_n$ involved in (3). According to Theorem 1 in the sequel, for consistency of the fence, one needs $c_n \to \infty$ at a certain rate, but there are many $c_n$'s that satisfy this requirement. Also note that although for the consistency it is not required that $\hat{\sigma}_{M,\tilde{M}}$ be a consistent estimator of $\sigma_{M,\tilde{M}}$ as long as it has the right order (see Section 4), there is always a constant involved which may make a difference in a finite sample situation. The problem can be solved by choosing a suitable $c_n$.

We now introduce the idea of an adaptive procedure. Recall that $\mathcal{M}$ denotes the set of candidate models, which includes a true model. To be more specific, we assume that there is a full model $M_f \in \mathcal{M}$, hence $\tilde{M} = M_f$ in (3); and that every model in $\mathcal{M} \setminus \{M_f\}$ is a submodel of a model in $\mathcal{M}$ with one less parameter than $M_f$. Let $M_*$ denote a model with minimum dimension among $M \in \mathcal{M}$. First note that ideally, one wishes to select $c_n$ that maximizes the probability of choosing the optimal model. Here for simplicity, the optimal model is defined as a true model that has the minimum dimension among all true models. This means that one wishes to choose $c_n$ that maximizes

$$(4) \qquad P = \mathrm{P}(M_0 = M_{\mathrm{opt}}),$$

where $M_{\mathrm{opt}}$ represents the optimal model and $M_0 = M_0(c_n)$ is the model selected by the fence procedure with the given $c_n$. However, two things are unknown on the right-hand side of (4): (i) under what distribution should the probability P be computed? and (ii) what is $M_{\mathrm{opt}}$?

To solve problem (i), note that the assumptions above on $\mathcal{M}$ imply that $M_f$ is a true model. Therefore, it is possible to bootstrap under $M_f$. For example, one may estimate the parameters under $M_f$, then use a model-based bootstrap to draw samples under $M_f$. This allows us to approximate the probability distribution P on the right side of (4).

To solve problem (ii), we use the idea of maximum likelihood. Namely, let $p^*(M) = \mathrm{P}^*(M_0 = M)$, where $M \in \mathcal{M}$ and $\mathrm{P}^*$ denotes the empirical probability obtained by bootstrapping. In other words, $p^*(M)$ is the sample proportion of times out of the total number of bootstrap samples that model $M$ is selected by the fence method with the given $c_n$. Let $p^* = \max_{M \in \mathcal{M}} p^*(M)$. Note that $p^*$ depends on $c_n$. The idea is to choose $c_n$ that maximizes $p^*$. It should be kept in mind that the maximization is not without restriction. To see this, note that if $c_n = 0$ then $p^* = 1$ (because when $c_n = 0$ the procedure always chooses $M_f$). Similarly, $p^* = 1$ for very large $c_n$, if $M_*$ is unique (because when $c_n$ is large enough the procedure always chooses $M_*$). Therefore, what one looks for is "the peak in the middle" of the plot of $p^*$ against $c_n$.

Here is another look at the method. Typically, the optimal model is the model from which the data is generated, then this model should be the most likely given the data. Thus, given $c_n$, one is looking for the model (using the fence procedure) that is most supported by the data or, in other words, one that has the highest (posterior) probability. The latter is estimated by bootstrapping. Note that although the bootstrap samples are generated under $M_f$, they are almost the same as those generated under the optimal model. This is because the estimates corresponding to the zero parameters are expected to be close to zero, provided that the parameter estimators under $M_f$ are consistent. (Note that in some special cases, a nonmodel based bootstrap algorithm can also be used. For instance, in the case of crossed random effects, Owen [22] presents a pigeonhole bootstrap algorithm which can be used effectively.) One then pulls off the $c_n$ that maximizes the (posterior) probability and this is the optimal choice, denoted by $c_n^*$.

A few technical issues deserve some attention:

1. Quite often the search for the peak in the middle finds multiple peaks (see Figure 1). In such cases, one should pick the highest. This is supported by our theoretical result, namely, Theorem 3 in the sequel which shows that $p^*(c_n^*) \to 1$ in probability as $n \to \infty$. It is also very common to have interval(s) of $c_n$ at which $p^*$ is at the maximum, say, $p^* = 1$. We then take the median of each interval, and let $c_n^*$ be the smallest of those medians, if there are more than one. The latest strategy is called *conservative*. For example, in the case of variable selection this strategy intends to make sure that no important variable is missing (in other words, to minimize the probability of underfit).

2. There are two extreme cases which occur when the optimal model is either the full model, $M_f$, or the minimum model, $M_*$. It should be pointed out that these cases are rare in practice. For example, in most cases of variable selection, there are a set of candidate variables and only some of them are important. This means that the optimal model is neither the full model nor the minimum model. Furthermore, when the extreme cases do occur, they are often easy to identify from the plot of $p^*$ (see Figure 1). Alternatively, one can run screen tests for the extreme cases. Such tests are recommended as supplementary tools to the inspection of the plot.

The first screen test is called full model test. The idea is the following. Define $\mathcal{M}_{f-1}$ as the set of all models with one less parameter than $M_f$.

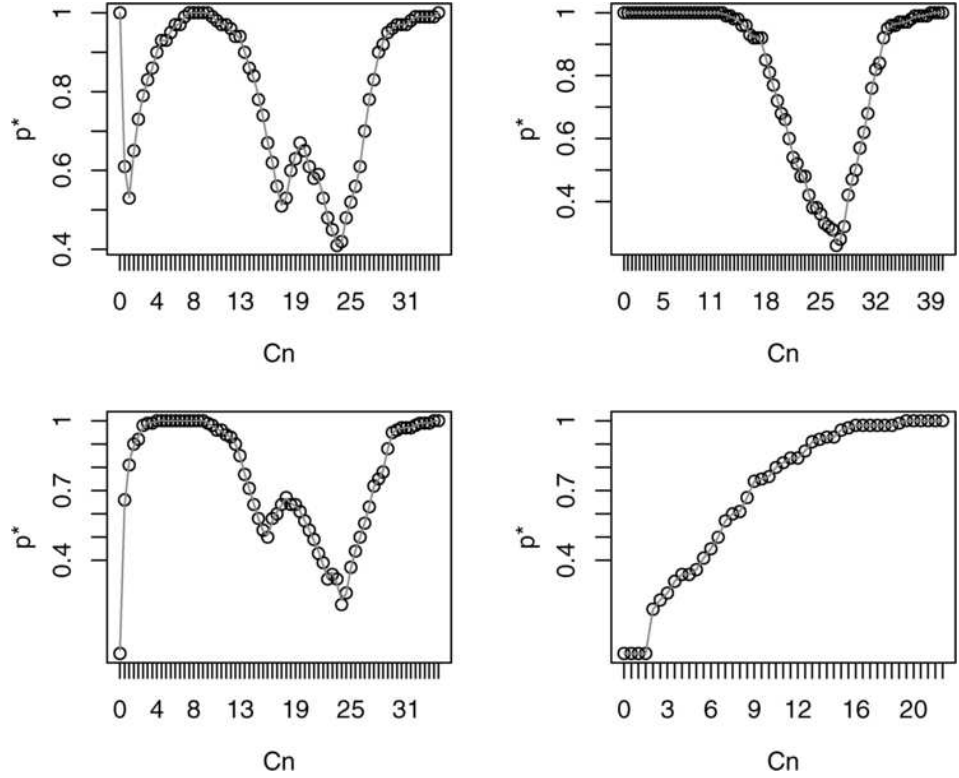

Fig. 1. *Upper left: $p^*$ versus $c_n$ without adjusting the baseline under Model 4. Upper right: $p^*$ versus $c_n$ without adjusting the baseline under Model 5. Lower left: $p^*$ versus $c_n$ with adjusting the baseline under Model 5. Lower right: $p^*$ versus $c_n$ with adjusting the baseline under Model 1. All plots are made based on 100 bootstrap samples generated given the first simulated dataset.*

Suppose that when $M_{\mathrm{f}}$ is the optimal model, we have $\mathrm{E}(\hat{Q}_M - \hat{Q}_{M_{\mathrm{f}}}) \sim a_n$, $\forall M \in \mathcal{M}_{\mathrm{f}-1}$. Here $u_n \sim v_n$ means that both $u_n/v_n$ and $v_n/u_n$ are bounded. On the other hand, if $M_{\mathrm{f}}$ is not optimal, there is $M \in \mathcal{M}_{\mathrm{f}-1}$ which is a true model, hence $\mathrm{E}(\hat{Q}_M - \hat{Q}_{M_{\mathrm{f}}}) = O(b_n)$, where $b_n = o(a_n)$. It follows that $\min_{M \in \mathcal{M}_{\mathrm{f}-1}} \mathrm{E}(\hat{Q}_M - \hat{Q}_{M_{\mathrm{f}}}) = O(b_n)$. Therefore, we consider

$$(5) \qquad q_n = \frac{\{\min_{M \in \mathcal{M}_{\mathrm{f}-1}} \mathrm{E}(\hat{Q}_M - \hat{Q}_{M_{\mathrm{f}}})\}^2}{a_n b_n}.$$

In practice, $q_n$ is replaced by its bootstrap estimate, $q_n^*$, obtained as above. If $q_n^* < 1$, the full model test passes; otherwise, the full model test fails, in which case we assign $c_n^* = 0$. The second screen test is called minimum model test. For simplicity, we assume that there is a unique $M_* \in \mathcal{M}$ that has the minimum dimension. Suppose that $\mathrm{E}(\hat{Q}_{M_*} - \hat{Q}_{M_{\mathrm{f}}}) = O(g_n)$ if $M_*$ is incorrect; and the order becomes $O(h_n)$ if $M_*$ is correct (hence optimal), where $h_n = o(g_n)$. We then consider

$$r_n = \frac{\{\text{E}(\hat{Q}_{M_*} - \hat{Q}_{M_\text{f}})\}^2}{g_n h_n}. \tag{6}$$

Let $r_n^*$ be the bootstrap version of $r_n$. If $r_n^* > 1$, the minimum model test passes; otherwise, the minimum model test fails, in which case we assign $c_n^*$ as the upper bound of a sequence of values considered (see below).

One concern about the screen tests is that the quantities $a_n$, $b_n$, $g_n$, $h_n$ may subject to scale change. Throughout this paper, we choose those quantities naturally without additional constants. For example, if $a_n = O(n)$, we simply take $a_n = n$ (not $2n$ or $3n$). On the other hand, the minimum model test can be replaced by the following threshhold checking which does not suffer from the scale change. Assuming that $M_*$ is true (therefore, optimal), one can draw bootstrap samples $y_{*,b}^*$, $b = 1, \ldots, B$ under $M_*$. Then based on such bootstrap samples, compute $d_* = \max_{1 \leq b \leq B} \{\hat{Q}_{M_*,b}^* - \hat{Q}_{M_\text{f}^*,b}^*\}$, where $M_\text{f}^*$ is defined below. If $\hat{Q}_{M_*} - \hat{Q}_{M_\text{f}^*} > d_*$, do not consider the right tail of the plot of $p^*$ against $c_n$ that goes up and stays at one (see Figure 1); otherwise, consider it. Unfortunately, the same idea does not apply to the full model case. To see why, note that the threshhold checking is similar to hypothesis testing. In the minimum model case, the null hypothesis is that $M_*$ is true, therefore, one can draw bootstrap samples under the null. However, in the full model case, the null hypothesis is that $M_\text{f}$ is optimal, which is equivalent to that none of the models in $\mathcal{M}_{f-1}$ [defined above (5)] is true. We do not know how to draw bootstrap samples under such a null. To solve this problem, we use a method called adjusting the baseline. Consider, for simplicity, the problem of selecting the fixed covariates under a linear mixed model. Suppose that the candidate variables are $X_1, \ldots, X_s$. Create an additional variable that is unrelated to the data, for example, by generating a random vector $X_{s+1}^*$ whose components are i.i.d. $\sim N(0,1)$ and are independent of the data. Define the model $M_\text{f}^*$ as the model that includes $X_1, \ldots, X_s, X_{s+1}^*$. Then replace $\hat{Q}_{M_\text{f}}$ in the fence inequality by $\hat{Q}_{M_\text{f}^*}$. Note that even though the baseline is adjusted, $M_\text{f}^*$ is *not* considered as a candidate model (because we know it is not optimal). Note that after the baseline change, $p^*$ will not equal to one when $c_n = 0$ (see Figure 1). Although the standard normal distribution is used to adjust the baseline, our simulation results (see Section 5.1) suggest that the method is quite stable with respect to different choices of baselines.

REMARK. In practice, if there is belief that $M_*$ and $M_\text{f}$ are unlikely to be the optimal model, neither the screen tests nor the baseline adjustment/threshhold checking are necessary.

3. Finally, one needs to determine at which values of $c_n$ to evaluate $p^*$. Theoretically, the range of $c_n$ is $[0, \infty)$, but practically one needs an upper bound. This can be determined as follows. Note that any $c_n \geq B = (\hat{Q}_{M_*} - \hat{Q}_{M_f})/\hat{\sigma}_{M_*, M_f}$ makes no difference to the fence procedure (assuming no baseline adjustment). This is because then (3) is satisfied by $M_*$, hence $M_0 = M_*$. Therefore, we choose the upper bound of $c_n$ as $B^* = [B] + 1$. We then divide the interval $[0, B^*]$ by subintervals of equal length and consider the end points.

REMARK. It turns out that requiring the existence of a full model or other known true model from which to draw bootstrap samples is not much of a practical problem, because in essence the adaptive fence can be done in two steps. In the first step, one could use the fence with a fixed $c_n$ (e.g., $c_n = 1$) to select a true model (which may not be optimal). Then in the second step, one applies the adaptive fence procedure with bootstrap samples drawn under the true model selected in the first step. Note that in the first step, one does not need $c_n$ to increase in order to select (with probability tending to one) a true model.

**4. Consistency of fence, F-B fence and adaptive fence.** We assume that the following Assumptions A1–A4 hold for each $M \in \mathcal{M}$, where $\theta_M$ represents a parameter vector at which $\mathrm{E}(Q_M)$ attains its minimum, and $\partial Q_M/\partial \theta_M$, and so forth, represent derivatives evaluated at $\theta_M$. Similarly, $\partial \tilde{Q}_M/\partial \theta_M$, and so forth, represent derivatives evaluated at $\tilde{\theta}_M$.

ASSUMPTION A1. $Q_M$ is three-times continuously differentiable with respect to $\theta_M$; and

$$\mathrm{E}\left(\frac{\partial Q_M}{\partial \theta_M}\right) = 0. \tag{7}$$

ASSUMPTION A2. There is a constant $B_M$ such that $Q_M(\tilde{\theta}_M) > Q_M(\theta_M)$, if $|\tilde{\theta}_M| > B_M$.

ASSUMPTION A3. The equation $\partial Q_M/\partial \theta_M = 0$ has an unique solution.

ASSUMPTION A4. There is a sequence of positive numbers $a_n \to \infty$ and $0 \leq \gamma < 1$ such that the following hold: $\partial Q_M/\partial \theta_M - \mathrm{E}(\partial Q_M/\partial \theta_M) = O_\mathrm{P}(a_n^\gamma)$, $\partial^2 Q_M/\partial \theta_M \, \partial \theta'_M - \mathrm{E}(\partial^2 Q_M/\partial \theta_M \, \partial \theta'_M) = O_\mathrm{P}(a_n^\gamma)$, $\liminf a_n^{-1} \lambda_{\min}\{\mathrm{E} \times (\partial^2 Q_M/\partial \theta_M \, \partial \theta'_M)\} > 0$, $\limsup a_n^{-1} \lambda_{\max}\{\mathrm{E}(\partial^2 Q_M/\partial \theta_M \, \partial \theta'_M)\} < \infty$, and there is $\delta_M > 0$ such that $\sup_{|\tilde{\theta}_M - \theta_M| \leq \delta_M} |\partial^3 \tilde{Q}_M/\partial \theta_{M,j} \, \partial \theta_{M,k} \, \partial \theta_{M,l}| = O_\mathrm{P}(a_n)$, $1 \leq j, k, l \leq p_M$, where $p_M = \dim(\theta_M)$.

In addition, we assume the following. Recall that $c_n$ is the constant in (3).

ASSUMPTION A5. $c_n \to \infty$; $\forall$ true model $M^*$ and incorrect model $M$, we have $\mathrm{E}(Q_M) > \mathrm{E}(Q_{M^*})$, $\liminf(\sigma_{M,M^*}/a_n^{2\gamma-1}) > 0$ and $c_n \sigma_{M,M^*}/\{\mathrm{E}(Q_M) - \mathrm{E}(Q_{M^*})\} \to 0$.

ASSUMPTION A6. $\hat{\sigma}_{M,M^*} > 0$ and $\hat{\sigma}_{M,M^*} = \sigma_{M,M^*} O_{\mathrm{P}}(1)$ if $M^*$ is true and $M$ incorrect; and $\sigma_{M,M^*} \vee a_n^{2\gamma-1} = \hat{\sigma}_{M,M^*} O_{\mathrm{P}}(1)$ if both $M$ and $M^*$ are true.

NOTE. (7) is satisfied if $\mathrm{E}(Q_M)$ can be differentiated inside the expectation. Assumption A2 implies that $|\hat{\theta}_M| \leq B_M$. To illustrate Assumptions A4 and A5, consider the case of clustered responses (see the last paragraph of Section 2.3). Then under regularity conditions, Assumption A4 holds with $a_n = m$ and $\gamma = 1/2$. Furthermore, we have $\sigma_{M,M^*} = O(\sqrt{m})$ and $\mathrm{E}(Q_M) - \mathrm{E}(Q_{M^*}) = O(m)$, provided that $M^*$ is true, $M$ is incorrect and some regularity conditions hold. Thus, Assumption A5 holds with $\gamma = 1/2$ and $c_n$ being any sequence satisfying $c_n \to \infty$ and $c_n/\sqrt{m} \to 0$. Finally, Assumption A6 does not require that $\hat{\sigma}_{M,M^*}$ be a consistent estimator of $\sigma_{M,M^*}$—only that it has the same order as $\sigma_{M,M^*}$.

LEMMA 2. *Under Assumptions A1–A4, we have $\hat{\theta}_M - \theta_M = O_{\mathrm{P}}(a_n^{\gamma-1})$ and $\hat{Q}_M - Q_M = O_{\mathrm{P}}(a_n^{2\gamma-1})$.*

Let $M_0$ be the model selected by fence using (3). The following theorem establishes consistency of the fence procedure.

THEOREM 1. *Under Assumptions A1–A6, we have with probability tending to one that $M_0$ is a true model with minimum dimension.*

The proofs of Lemma 2 and Theorem 1 are given in Sections 7.1 and 7.2, respectively.

The next theorem establishes consistency of the F-B fence proposed in Section 3.1. Note that the method is introduced in the case of extended GLMMs. Let $M_0^{\dagger}$ be the final model selected by the F-B fence procedure using (3).

THEOREM 2. *Under Assumptions A1–A6, we have with probability tending to one that $M_0^{\dagger}$ is a true model and no proper submodel of $M_0^{\dagger}$ is a true model.*

Note that here consistency is in the sense that with probability tending to one, $M_0^\dagger$ is a true model which cannot be further reduced or simplified. The proof is given in Section 7.3.

Finally, we give sufficient conditions for the consistency of the adaptive fence procedure introduced in Section 3.2. For simplicity, assume that $M_{\text{opt}}$ is unique. Consider the ratios $r_M = (\hat{Q}_M - \hat{Q}_{M_{\text{f}}})/\hat{\sigma}_{M,M_{\text{f}}}$, $M \in \mathcal{M}$. Let $\mathcal{M}_{\text{w}\leq}$ denote the subset of incorrect models with dimension $\leq |M_{\text{opt}}|$. Write $r_{\text{opt}} = r_{M_{\text{opt}}}$ and $r_{\text{w}\leq} = \min_{M \in \mathcal{M}_{\text{w}\leq}} r_M$. Denote the c.d.f.s of $r_{\text{opt}}$ and $r_{\text{w}\leq}$ by $F_{\text{opt}}$ and $F_{\text{w}\leq}$, respectively. Let $M_0(x)$ be the model selected by the fence procedure using (3) with $c_n = x$, and $P(x) = \mathrm{P}(M_0(x) = M_{\text{opt}})$. Let $P^*(x)$ be the bootstrap version of $P(x)$. Denote the bootstrap sample size by $n^*$. Recall the definitions of $a_n$, $b_n$, $q_n$, $q_n^*$ in (5), $g_n$, $h_n$, $r_n$, $r_n^*$ in (6), and $B^*$ above the final remark of Section 3.2. We make the following assumptions.

ASSUMPTION A7 (*Asymptotic distributional separation*). If $M_{\text{opt}} \notin \{M_{\text{f}}, M_*\}$, then for any $\varepsilon > 0$, there is $0 < \delta \leq 0.1$, $x_{n,1} < x_{n,2} < x_{n,3}$, and $N \geq 1$ such that when $n \geq N$ the following hold: $F_{\text{opt}}(x_{n,1}) > 1 - \varepsilon$, $F_{\text{w}\leq}(x_{n,3}) \leq \varepsilon$, $P(x_{n,2}) > 1 - \delta$, $1 - 4\delta < P(x_{n,j}) \leq 1 - 3\delta$, $j = 1, 3$; if $M_{\text{opt}} = M_{\text{f}}$, we have $\mathrm{P}(\min_{M \in \mathcal{M}, M \neq M_{\text{f}}} \hat{Q}_M > \hat{Q}_{M_{\text{f}}}) \to 1$ as $n \to \infty$.

ASSUMPTION A8 (*Good bootstrap approximation*). If $M_{\text{opt}} \notin \{M_{\text{f}}, M_*\}$, then for any $\delta, \eta > 0$, there are $N \geq 1$, $N^* = N^*(n)$ such that, when $n \geq N$ and $n^* \geq N^*$, we have $\mathrm{P}(\sup_{x>0} |P^*(x) - P(x)| < \delta) > 1 - \eta$; if $M_{\text{opt}} = M_{\text{f}}$, we have $q_n/q_n^* = O_{\mathrm{P}}(1)$; if $M_{\text{opt}} = M_*$, we have $q_n^*/q_n = O_{\mathrm{P}}(1)$ and $r_n^*/r_n = O_{\mathrm{P}}(1)$.

For the most part, Assumption A7 says that there is an asymptotic separation between the optimal model and the incorrect ones that matter in that the peak of $P(x)$ is distant from the area where $r_{\text{w}\leq}$ concentrates. This is reasonable because, typically, $r_{\text{opt}}$ is of lower order than $r_{\text{w}\leq}$. Therefore, one can find an interval, $(x_{n,1}, x_{n,3})$, such that (3) is almost always satisfied by $M = M_{\text{opt}}$ when $c_n \in (x_{n,1}, x_{n,3})$. On the other hand, $(x_{n,1}, x_{n,3})$ is distant from the area where $r_{\text{w}\leq}$ concentrates, so that $r_{\text{opt}} \leq c_n$, $r_{\text{w}\leq} > c_n$ with high probability, if $c_n \in (x_{n,1}, x_{n,3})$. Thus, $P(x)$ is expected to peak in $(x_{n,1}, x_{n,3})$ while $F_{\text{w}\leq}(x)$ stays low in the region.

Recall that $p^*$ in the adaptive procedure is a function of $c_n$, that is, $p^* = p^*(c_n)$. The following theorem establishes consistency of the adaptive fence. The proof is given in Section 7.4.

THEOREM 3. *Under Assumptions* A7 *and* A8, *the following hold:*

(i) *If* $M_{\text{opt}} \notin \{M_{\text{f}}, M_*\}$, *then with probability tending to one there is* $c_n^* \in (0, \infty)$ *which is at least a local maximum and approximate global maximum* of $p^*$ in the sense that for any $\delta, \eta > 0$, there is $N \geq 1$ and $N^* = N^*(n)$ such that $\mathrm{P}(p^*(c_n^*) \geq 1 - \delta) \geq 1 - \eta$, if $n \geq N$ and $n^* \geq N^*$.

(ii) *In general, define $c_n^*$ as*

$$\begin{cases} 0, & \text{if } q_n^* > 1, \\ B^*, & \text{if } q_n^* \leq 1, r_n^* < 1, \\ \text{the } c_n^* \text{ in } (\mathrm{i}), & \text{if } q_n^* \leq 1, r_n^* \geq 1 \text{ and such a } c_n^* \text{ exists}, \\ 1, & \text{otherwise}. \end{cases}$$

Let $M_0^*$ be the model selected by the fence procedure using (3) with $\tilde{M} = M_\mathrm{f}$ and $c_n$ replaced by $c_n^*$. Then $M_0^*$ is consistent in the sense that for any $\eta > 0$ there is $N \geq 1$, $N^* = N^*(n)$ such that $\mathrm{P}(M_0^* = M_\mathrm{opt}) \geq 1 - \eta$, if $n \geq N$ and $n^* \geq N^*$.

## 5. Examples of simulations and data analysis.

5.1. *The Fay–Herriot model—an illustration of adaptive fence method.* The Fay–Herriot model is widely used in small area estimation. It was first proposed to estimate the per-capita income of small places with population less than 1,000 (Fay and Herriot [8]). The model can expressed as $y_i = x_i'\beta + v_i + e_i$, $i = 1, \ldots, m$, where $x_i$ is a vector of known covariates, $\beta$ is a vector of unknown regression coefficients, $v_i$'s are area-specific random effects and $e_i$'s represent sampling errors. It is assumed that $v_i$, $e_i$ are independent with $v_i \sim N(0, A)$ and $e_i \sim N(0, D_i)$. The variance $A$ is unknown, but the sampling variances $D_i$'s are assumed known.

Let $X = (x_i')_{1 \leq i \leq m}$, so that the model can be expressed as $y = X\beta + v + e$, where $y = (y_i)_{1 \leq i \leq m}$, $v = (v_i)_{1 \leq i \leq m}$ and $e = (e_i)_{1 \leq i \leq m}$. The first column of $X$ is assumed to be $1_m$ which corresponds to the intercept. The rest of the columns of $X$ are to be selected from a set of candidate covariate vectors $X_2, \ldots, X_K$, which include the true covariate vectors. First note that by applying the following transformation, we can simplify the problem to the case $D_i = 1$. Let $D = 1 + \max_{1 \leq i \leq m} D_i$. Draw independent samples $u_1, \ldots, u_m$ independent with the $v_i$'s and $e_i$'s such that $u_i \sim N(0, D - D_i)$, $1 \leq i \leq m$. Then let $\tilde{y}_i = (y_i + u_i)/\sqrt{D}$, $\tilde{x}_i = x_i/\sqrt{D}$, $\tilde{v}_i = v_i/\sqrt{D}$ and $\tilde{e}_i = (e_i + u_i)/\sqrt{D}$. Consider $\tilde{y}_i$'s as the new observations. Then we have $\tilde{y}_i = \tilde{x}_i'\beta + \tilde{v}_i + \tilde{e}_i$, $i = 1, \ldots, m$, where $\tilde{v}_i$, $\tilde{e}_i$, $i = 1, \ldots, m$, are independent with $\tilde{v}_i \sim N(0, \tilde{A})$, $\tilde{A} = A/D$ and $\tilde{e}_i \sim N(0, 1)$. Thus, without loss of generality, we let $D_i = 1$, $1 \leq i \leq m$.

Consider the fence ML model selection (see Section 2). It is easy to show that in this case, $\hat{Q}_M = (m/2)\{1 + \log(2\pi) + \log(|P_{X^\perp}y|^2/m)\}$, where $P_{X^\perp} = I_m - P_X$ and $P_X = X(X'X)^{-1}X'$. We assume for simplicity that $X$ is of full rank. Then $\hat{Q}_M - \hat{Q}_{M_\mathrm{f}} = (m/2)\log(|P_{X^\perp}y|^2/|P_{X_\mathrm{f}^\perp}y|^2)$. Furthermore, it can be shown that when $M$ is a true model, we have $\hat{Q}_M - \hat{Q}_{M_\mathrm{f}} =$

TABLE 1
*Fence methods with different $c_n$'s in the Fay–Herriot model*

| Optimal model | 1 | 2 | 3 | 4 | 5 |
|---|---|---|---|---|---|
| Adaptive $c_n$ (ST) | 100 | 100 | 100 | 99 | 100 |
| Adaptive $c_n$ (B/T) | 99 | 100 | 100 | 99 | 100 |
| $c_n = \log\log(n)$ | 52 | 63 | 70 | 83 | 100 |
| $c_n = \log(n)$ | 96 | 98 | 99 | 96 | 100 |
| $c_n = \sqrt{n}$ | 100 | 100 | 100 | 100 | 100 |
| $c_n = n/\log(n)$ | 100 | 91 | 95 | 90 | 100 |
| $c_n = n/\log\log(n)$ | 100 | 0 | 0 | 0 | 6 |

$(m/2)\log(1 + \frac{K-p}{m-K-1}F)$, where $p+1$ is the number of columns of $X$, and $F \sim F_{K-p,m-K-1}$. Therefore, $\sigma_{M,M_{\mathrm{f}}}$ is completely known given $|M|$ and can be evaluated accurately (e.g., by numerical integration).

We carry out a simulation study to evaluate the performance of the adaptive method. Here we consider the adaptive method assisted either by the screen tests (ST) or the baseline adjustment/threshhold checking (B/T). We consider a (relatively) small sample situation with $m = 30$. With $K = 5$, $X_2, \ldots, X_5$ were generated from the $N(0,1)$ distribution, and then fixed throughout the simulations. The candidate models include all possible models with at least an intercept (thus, there are $2^4 = 16$ candidate models). We consider five cases in which the data $y$ is generated from the model $y = \sum_{j=1}^{5} \beta_j X_j + v + e$, where $\beta' = (\beta_1, \ldots, \beta_5) = (1,0,0,0,0), (1,2,0,0,0), (1,2,3,0,0), (1,2,3,2,0)$ and $(1,2,3,2,3)$, denoted by Models 1, 2, 3, 4, 5, respectively. The true value of $A$ is 1 in all cases. The number of bootstrap samples for the evaluation of the $p^*$'s is 100.

In addition to the adaptive method, we consider five different (nonadaptive) $c_n$'s ($n = m$ in this case), which satisfy the consistency requirements given in Theorem 1 (note that these requirements reduce to $c_n \to \infty$ and $c_n/n \to 0$ in this case). These are $c_n = \log\log(n), \log(n), \sqrt{n}, n/\log(n)$ and $n/\log\log(n)$. Reported in Table 1 are percentage of times, out of 100 simulations that the optimal model was selected by each method.

*Summary.* Although the reported results for Adaptive $c_n$ (B/T) were obtained using $N(0,1)$ for the baseline adjustment, the same simulations were carried out when $N(0,1)$ is replaced by Uniform$[0,1]$, Poisson$(1)$ and Bernoulli distributions. The only (slight) differences in the results are those under Model 1, which are 99, 98 and 100, respectively, for Uniform$[0,1]$, Poisson$(1)$ and Bernoulli. This suggests that the method is not very sensitive to different choices of baselines which is what one desires. Figure 1 displays the plots of $p^*$ against $c_n$ in a number of situations. Furthermore, we explore the two-step adaptive fence procedure (with ST) described in the last remark of Section 3.2 and the same results were obtained.

It seems that performance of the fence with $c_n = \log(n)$, $\sqrt{n}$ or $n/\log(n)$ is fairly close to that of the adaptive fence. In any particular situation, one might get lucky to find a good $c_n$ value by chance, but one cannot be lucky all the time. In fact, for more complicated mixed models, the definition of the sample size may not simply be the total number of observations or the number of clusters so, for example, something like $\log(n)$ or $\sqrt{n}$ may not make sense.

COMPUTATIONAL NOTE. The simulations of this subsection were run on a Pentium Dual Core CPU 3.2 GHz, memory 4 GB, Harddrive 500 GB. The times it took to run the first simulation of Adaptive $c_n$ (B/T) under Models 1–5 were 1.7 sec., 3.0 sec., 4.1 sec., 4.4 sec. and 5.3 sec., respectively.

5.2. *Linear mixed models for clustered data.* We consider the following linear mixed model (see Jiang and Rao [16]), $y_{ij} = x'_{ij}\beta + \alpha_i + \varepsilon_{ij}$, $i = 1, \ldots, m$, $j = 1, \ldots, K$, where $x_{ij}$ is a vector of covariates and $\beta$ a vector of unknown regression coefficients (the fixed effects). The random effects $\alpha_1, \ldots, \alpha_m$, are generated independently from $N(0, \sigma^2)$. The errors are generated so that $\varepsilon_i = (\varepsilon_{ij})_{1 \leq j \leq K}$, $i = 1, \ldots, m$, are independent and multivariate normal with $\text{Var}(\varepsilon_i) = \tau^2\{(1-\rho)I + \rho J\}$, where $I$ is the identity matrix and $J$ matrix of 1's. Finally, the random effects are uncorrelated with the errors.

Now pretend that the covariance matrix of the data is unknown. The problem is to select the fixed covariates. Write the model as $y = X\beta + Z\alpha + \varepsilon$. The candidate covariates which are columns of $X$ are $X_1, \ldots, X_5$, where $X_1$ is a vector of 1's and $X_2, \ldots, X_4$ are generated randomly from the $N(0, 1)$ distribution, and then fixed throughout the simulations. We consider the $Q_M$ for LS model selection (described above Section 2.1) which is suitable for this situation.

We examine the performance of fence with fixed $c_n = 1.1$ and that of the adaptive fence. As comparison, two GICs developed in Jiang and Rao [16]) are considered, which are similar to (1) for this problem: (i) $\lambda_n = 2$ which corresponds to the $C_p$ method; and (ii) $\lambda_n = \log n$, where $n = mK$, which corresponds to the BIC. The latter choice satisfies the conditions of Theorem 1 in Jiang and Rao [16] for consistent model selection for this setting.

We consider the case where the errors have varying degrees of exchangeable structure. Four values of $\rho$ were considered: $0, 0.2, 0.5, 0.8$. The random effects and errors were simulated from normal distributions with $\sigma = \tau = 1$. We set the number of clusters to be $m = 100$ and the number of observations within a cluster to be $K = 5$. Three (true) $\beta$'s are considered: $(2, 0, 0, 4, 0)$, $(2, 9, 0, 4, 8)$ and $(1, 2, 3, 2, 3)$. A total of 100 realizations of each simulation were run.

*Summary.* The results reported in Table 2 for adaptive fence are those under B/T (see the previous subsection). The same results were obtained under ST. The fence method with fixed $c_n$ is seen to have robust selection performance in most situations considered. In cases where the true model was relatively small in dimension, the fence method suffers some from overfitting. The overfitting proneness in these few situations is less than that found when using $C_p$ but more than that found when using BIC. Selection performance in the second situation with a larger true model with high signal is solid for the fence method. However, in the last situation with the optimal model being the full model with all weak covariates, both BIC and $C_p$ tend to underfit. The fence method still shines having excellent performance with comparatively little or no underfitting empirically observed (note that overfitting is not possible in this situation). The effect of increasing correlation in the errors (i.e., clustering) is to act as a means of reducing effective sample size. The end result is that as the correlation between observations within a cluster increases, selection performance for all fixed penalization methods degrades somewhat. The adaptive fence on the other hand shines in all situations giving 100% selection accuracy. This clearly demonstrates the effectiveness of the adaptive fence method (at a computational cost, of course).

5.3. *Prenatal care for pregnancy.* This real-data example is an application of the F-B fence procedure to GLMMs (see Section 3.1). Rodriguez and Goldman [25] considered a dataset from a survey conducted in Guatemala regarding the use of modern prenatal care for pregnancies where some form of care was used (Pebley [23]). While Rodriguez and Goldman focused on assessing the performance of the approximation method that they developed in fitting a three-level variance component logistic model, we consider applying the fence method in selection of the fixed covariates in the variance component logistic model. The models are described as follows.

TABLE 2
*Simulation results: linear mixed model selection. Reported are probabilities of correct selection (underfitting, overfitting) as percentages estimated empirically from 100 realizations of the simulation. $C_p$ and BIC results for Models 1 and 2 were taken from Jiang and Rao [16]*

| Optimal model | $\rho$ | $C_p$ | BIC | Fence ($c_n = 1.1$) | Adaptive fence |
|---|---|---|---|---|---|
| $\beta' = (2,0,0,4,0)$ | 0 | 64(0, 36) | 97(0, 3) | 94(0, 6) | 100(0, 0) |
| | 0.2 | 57(0, 43) | 94(0, 6) | 91(0, 9) | 100(0, 0) |
| | 0.5 | 58(0, 42) | 96(1, 3) | 86(0, 14) | 100(0, 0) |
| | 0.8 | 61(0, 39) | 96(0, 4) | 72(0, 28) | 100(0, 0) |
| $\beta' = (2,9,0,4,8)$ | 0 | 87(0, 13) | 99(0, 1) | 100(0, 0) | 100(0, 0) |
| | 0.2 | 87(0, 13) | 99(0, 1) | 100(0, 0) | 100(0, 0) |
| | 0.5 | 80(0, 20) | 99(0, 1) | 99(0, 1) | 100(0, 0) |
| | 0.8 | 78(1, 21) | 96(1, 3) | 94(0, 6) | 100(0, 0) |
| $\beta' = (1,2,3,2,3)$ | 0 | 85(15, 0) | 81(19, 0) | 100(0, 0) | 100(0, 0) |
| | 0.2 | 79(21, 0) | 73(27, 0) | 100(0, 0) | 100(0, 0) |
| | 0.5 | 74(26, 0) | 64(36, 0) | 97(3, 0) | 100(0, 0) |
| | 0.8 | 44(56, 0) | 26(74, 0) | 94(6, 0) | 100(0, 0) |

Suppose that given the random effects at community levels $u_i$, $1 \leq i \leq m$ and random effects at family levels $v_{ij}$, $1 \leq i \leq m$, $1 \leq j \leq n_i$, binary responses $y_{ijk}$, $1 \leq i \leq m$, $1 \leq j \leq n_i$, $1 \leq k \leq n_{ij}$, are conditionally independent with $\pi_{ijk} = \mathrm{E}(y_{ijk}|u,v) = \mathrm{P}(y_{ijk}=1|u,v)$. Furthermore, suppose that the random effects are independent with $u_i \sim N(0,\sigma^2)$ and $v_{ij} \sim N(0,\tau^2)$. The following models for the conditional means are considered such that under model $M$, $\mathrm{logit}(\pi_{ijk}) = X'_{M,ijk}\beta_M + u_i + v_{ij}$, where $X_{M,ijk}$ is a subvector of the full set of fixed covariates and $\beta_M$ the corresponding vector of regression coefficients.

Let $\psi = (\sigma^2, \tau^2)'$. The vector of parameters under model $M$ is $\theta_M = (\beta'_M, \psi')'$. We use the $Q_M$ introduced earlier for extended GLMMs (see the second paragraph of Section 3.1). An estimated $\sigma^2_{M,M^*}$ can be obtained using the idea of observed variance (see Section 2.3, and Jiang et al. [17] for detail). The expectations involved in $Q_M$ are evaluated by numerical integration. Since the number of covariates considered is quite large, to keep the computational time manageable, we apply the F-B fence procedure introduced in Section 3.1 with $c_n = 1$.

The data analysis has selected the following variables (in the order that they were selected in the forward procedure): Proportion indigenous (1981), Modern toilet in household, Husband's education secondary or better, Husband's education primary, Television watched daily, Distance to nearest clinic, Mother's education primary, Television not watched daily, Mother's education secondary or better, Indigenous (no Spanish), Indigenous (Spanish), Mother age, Husband agriculture employee, Husband agriculture self-employee, Child age, Birth order 4–6 and Husband's education missing. There are some interesting differences between the fixed effects discovered by the fence versus those found by standard maximum likelihood analysis using a 5% significance level as reported in Rodriguez and Goldman [25]. First, Husband's education overall (primary or higher relative to the reference group of no education for the husband) was found to be an important predictor whereas Rodriguez and Goldman found that only Husband's secondary education was important. Our more uniform finding is also in line with the finding for Mother's education. The implication is that education of some kind is important for both the mother and husband to have. A similar kind of finding was observed for variables corresponding to husband's profession. We found that regardless of what type of agricultural employment the husband had, it was an important predictor overall. Rodriguez and Goldman report that only nonself employed agricultural jobs for the husband mattered. The fence method also uniquely found that watching television (daily or not) was an important predictor. This can be intuitively justified since it provides a medium for women to learn more about modern prenatal health care methods, and thus make it more likely for them to choose to use such methods. Other findings were in line with those of Rodriguez and Goldman [25].

**6. Concluding remarks.** Fence is different from procedures like AIC, BIC in that there is no criterion function that is minimized. In other words, instead of trying to find an "optimal" model that minimizes a criterion function, fence proposes to carry out the optimization by two steps. The first step is to identify the set of true models (the ones that are in the fence) or, in case a true model does not exist, the models that best approximate the real-life problem. Note that although in this paper we have assumed the existence of a true model, the method can be easily extended to the situation where a true model does not exist, or is understood as the one that provides the best approximation. On the other hand, the second step of fence, which identifies the model with minimal dimension within the fence, is quite flexible. For example, the dimension of a model may not be defined as the number of estimated parameters (e.g., Hastie and Tibshirani [12], Ye [29]); or it may be replaced by some other considerations, such as economical concerns. In fact, practically speaking, optimality in model selection usually goes beyond statistics. Keeping this in mind, it appears that the fence procedure is easier to incorporate with other scientific or economical criteria than minimizing a single criterion function determined before the scientific or economic problem.

A good feature of the fence algorithm is that it needs not search over all the candidate models in order to find the optimal model.

In this paper, we have demonstrated the robust performance of fence in linear and generalized linear mixed model selection. In addition, we have introduced a stepwise fence procedure to handle situations of large number of predictors. Furthermore, we have proposed an adaptive procedure for choosing a tuning constant involved in the fence method. The adaptive procedure improves the finite sample performance of fence at a computational cost. On the theoretical side, we have established consistency of the different fence procedures, with the proofs given in the next section.

## 7. Proofs.

7.1. *Proof of Lemma* 2. Assumptions A2 and A3 imply that $\hat{\theta}_M$ is the unique solution to $\partial Q_M / \partial \theta_M = 0$. By Taylor expansion, we have

$$\tilde{Q}_M - Q_M$$
$$= \left(\frac{\partial Q_M}{\partial \theta_M}\right)'(\tilde{\theta}_M - \theta_M) + \frac{1}{2}(\tilde{\theta}_M - \theta_M)'\left(\frac{\partial^2 Q_M}{\partial \theta_M \partial \theta_M'}\right)(\tilde{\theta}_M - \theta_M)$$
$$+ \frac{1}{6}\sum_{j,k,l}\left(\frac{\partial^3 Q_M^*}{\partial \theta_{M,j} \partial \theta_{M,k} \partial \theta_{M,l}}\right)(\tilde{\theta}_{M,j} - \theta_{M,j})(\tilde{\theta}_{M,k} - \theta_{M,k})(\tilde{\theta}_{M,l} - \theta_{M,l})$$
$$= I_1 + \frac{1}{2}I_2 + \frac{1}{6}I_3$$

for any $\tilde{\theta}_M$, where $\partial^3 Q_M^* / \cdots$ represents the third derivatives evaluated at $\theta_M^*$, which lies between $\theta_M$ and $\tilde{\theta}_M$. For any $\varepsilon > 0$, by Assumptions A1 and A4, there are $\delta > 0$ and $N_0 \geq 1$ such that $\lambda_{\min}\{\mathrm{E}(\partial^2 Q_M/\partial \theta_M \partial \theta_M')\} \geq \delta a_n$, $n \geq N_0$, and $L_1 > 0$ such that the probability is greater than $1 - \varepsilon$ that $|\partial Q_M/\partial \theta_M| \leq L_1 a_n^{\gamma}$,

$$\left\|\frac{\partial^2 Q_M}{\partial \theta_M \partial \theta_M'} - \mathrm{E}\left(\frac{\partial^2 Q_M}{\partial \theta_M \partial \theta_M'}\right)\right\| \leq L_1 a_n^{\gamma},$$

$$\max_{j,k,l} \sup_{|\tilde{\theta}_M - \theta_M| \leq \delta_M} \left|\frac{\partial^3 \tilde{Q}_M}{\partial \theta_{M,j} \partial \theta_{M,k} \partial \theta_{M,l}}\right| \leq L_1 a_n.$$

Now choose $L_2 > 0$ such that $\delta L_2 > 2L_1$. Let $\Theta_{M,L_2} = \{\tilde{\theta}_M : |\tilde{\theta}_M - \theta_M| \leq L_2 a_n^{\gamma-1}\}$, and $\bar{\Theta}_{M,L_2}$ be the boundary of $\Theta_{M,L_2}$, that is, $\bar{\Theta}_{M,L_2} = \{\tilde{\theta}_M : |\tilde{\theta}_M - \theta_M| = L_2 a_n^{\gamma-1}\}$. Then choose $N_1 \geq 1$ such that $L_2 a_n^{\gamma-1} \leq \delta_M$, $n \geq N_1$. It follows that for $\tilde{\theta} \in \bar{\Theta}_{M,L_2}$, we have $|I_1| \leq L_1 L_2 a_n^{2\gamma-1}$, $I_2 \geq \delta L_2^2 a_n^{2\gamma-1} - L_1 L_2^2 a_n^{3\gamma-2}$, $|I_3| \leq L_1 a_n (\sum_j |\tilde{\theta}_{M,j} - \theta_{M,j}|)^3 \leq L_1 L_2^3 p_M^{3/2} a_n^{3\gamma-2}$, hence

$$\tilde{Q}_M - Q_M$$
(8) $$\geq \frac{1}{2}L_2 a_n^{2\gamma-1}\{\delta L_2 - 2L_1 - L_1 L_2(1 + \frac{1}{3}L_2 p_M^{3/2})a_n^{\gamma-1}\},$$

$\forall \tilde{\theta} \in \bar{\Theta}_{M,L_2}$. If we choose $N_2 \geq 1$ such that, when $n \geq N_2$, the quantity inside $\{\cdots\}$ on the right-hand side of (8) is positive, and let $N = N_0 \vee N_1 \vee N_2$, then we have with probability greater than $1 - \varepsilon$, that $\tilde{Q}_M > Q_M$, $\forall \tilde{\theta} \in \bar{\Theta}_{M,L_2}$. It follows that $\mathrm{P}(|\hat{\theta}_M - \theta_M| < L_2 a_n^{\gamma-1}) \geq 1 - \varepsilon$, if $n \geq N$. This proves that $\hat{\theta}_M - \theta_M = O_\mathrm{P}(a_n^{\gamma-1})$.

By similar arguments, it can be shown that for any $\varepsilon > 0$, there are constants $L$, $L_1$, $L_2$ and $N \geq 1$ such that, when $n \geq N$,

$$\hat{Q}_M - Q_M \leq L_1 L_2 a_n^{2\gamma-1} + \frac{1}{2}LL_2^2 a_n^{2\gamma-1} + \frac{1}{2}L_1 L_2^2 a_n^{3\gamma-2} + \frac{1}{6}L_1 L_2^3 p_M^{3/2} a_n^{3\gamma-2}$$

$$\leq L_2\{L_1 + \tfrac{1}{2}(L+L_1)L_2 + \tfrac{1}{6}L_1L_2^2 p_M^{3/2}\}a_n^{2\gamma-1}$$

with probability $> 1 - \varepsilon$. This proves that $\hat{Q}_M - Q_M = O_{\mathrm{P}}(a_n^{2\gamma-1})$.

7.2. *Proof of Theorem* 1. For the most part, we show that with probability tending to one (w.p. $\to 1$), all the true models (with $|M| < |\tilde{M}|$) are in the fence, and all the incorrect ones are out.

Let $M$ be an incorrect model and $M^*$ a true model. By Lemma 2 and Assumption A5, we have $\hat{Q}_M - \hat{Q}_{M^*} = Q_M - Q_{M^*} + \hat{Q}_M - Q_M - (\hat{Q}_{M^*} - Q_{M^*}) = Q_M - Q_{M^*} + O_{\mathrm{P}}(a_n^{2\gamma-1}) = \mathrm{E}(Q_M) - \mathrm{E}(Q_{M^*}) + \{Q_M - Q_{M^*} - \mathrm{E}(Q_M - Q_{M^*})\} + O_{\mathrm{P}}(a_n^{2\gamma-1}) = \mathrm{E}(Q_M) - \mathrm{E}(Q_{M^*}) + \sigma_{M,M^*} O_{\mathrm{P}}(1) = \{\mathrm{E}(Q_M) - \mathrm{E}(Q_{M^*})\}\{1 + o_{\mathrm{P}}(1)\}$. It follows that, w.p. $\to 1$, we have $\hat{Q}_M > \hat{Q}_{M^*}$. This implies that w.p. $\to 1$, $\tilde{M}$ is a true model (because an incorrect model cannot be the minimizer).

Furthermore, it is seen from this argument that if $M$ is incorrect, we have

$$\hat{Q}_M - \hat{Q}_{M^*}$$
$$(9) \qquad = c_n \hat{\sigma}_{M,M^*} \left[ \frac{c_n \sigma_{M,M^*}}{\mathrm{E}(Q_M) - \mathrm{E}(Q_{M^*})} \left(\frac{\hat{\sigma}_{M,M^*}}{\sigma_{M,M^*}}\right)\{1 + o_{\mathrm{P}}(1)\}^{-1} \right]^{-1}.$$

Assumptions A5 and A6 imply that the quantity inside $[\cdots]$ in (9) is $o_{\mathrm{P}}(1)$. Therefore, w.p. $\to 1$, we have $\hat{Q}_M > \hat{Q}_{M^*} + c_n \hat{\sigma}_{M,M^*}$. It follows that

$$\mathrm{P}(|M| < |\tilde{M}|, M \in \tilde{\mathcal{M}}_-)$$
$$\leq \mathrm{P}(\hat{Q}_M \leq \hat{Q}_{\tilde{M}} + c_n \hat{\sigma}_{M,\tilde{M}})$$
$$\leq \sum_{M^* \text{ is true}} \mathrm{P}(\hat{Q}_M \leq \hat{Q}_{M^*} + c_n \hat{\sigma}_{M,M^*}, \tilde{M} = M^*) + \mathrm{P}(\tilde{M} \text{ is incorrect})$$
$$\leq \sum_{M^* \text{ is true}} \mathrm{P}(\hat{Q}_M \leq \hat{Q}_{M^*} + c_n \hat{\sigma}_{M,M^*}) + \mathrm{P}(\tilde{M} \text{ is incorrect}) \to 0.$$

Let $E_1 = \bigcap_{M \text{ is incorrect}, |M| < |\tilde{M}|} \{M \notin \tilde{\mathcal{M}}_-\}$, then $E_1^c = \bigcup_{M \text{ is incorrect}} \{|M| < |\tilde{M}|, M \in \tilde{\mathcal{M}}_-\}$, hence $\mathrm{P}(E_1^c) \to 0$. This proves the "out" part.

On the other hand, if $M$ and $M^*$ are both true models, then by the property of $Q_M$, we have $\mathrm{E}(Q_M) = \mathrm{E}(Q_{M^*})$. Therefore, by similar arguments and Assumption A6, we have $\hat{Q}_M - \hat{Q}_{M^*} = Q_M - Q_{M^*} + O_{\mathrm{P}}(a_n^{2\gamma-1}) = \hat{\sigma}_{M,M^*} O_{\mathrm{P}}(1)$. Since $c_n \to \infty$, we have, w.p. $\to 1$, $\hat{Q}_M \leq \hat{Q}_{M^*} + c_n \hat{\sigma}_{M,M^*}$. It follows that

$$\mathrm{P}(|M| < |\tilde{M}|, M \notin \tilde{\mathcal{M}}_-)$$
$$\leq \mathrm{P}(\hat{Q}_M > \hat{Q}_{\tilde{M}} + c_n \hat{\sigma}_{M,\tilde{M}})$$
$$\leq \sum_{M^* \text{ is true}} \mathrm{P}(\hat{Q}_M > \hat{Q}_{M^*} + c_n \hat{\sigma}_{M,M^*}, \tilde{M} = M^*) + \mathrm{P}(\tilde{M} \text{ is incorrect})$$

$$\leq \sum_{M^* \text{ is true}} \mathrm{P}(\hat{Q}_M > \hat{Q}_{M^*} + c_n \hat{\sigma}_{M,M^*}) + \mathrm{P}(\tilde{M} \text{ is incorrect}) \to 0.$$

Let $E_2 = \bigcap_{M \text{ is true}, |M| < |\tilde{M}|} \{M \in \tilde{\mathcal{M}}_-\}$, then $E_2^c = \bigcup_{M \text{ is true}} \{|M| < |\tilde{M}|, M \notin \tilde{\mathcal{M}}_-\}$, hence $\mathrm{P}(E_2^c) \to 0$. This proves the "in" part.

Finally, note that $\{M_0 \text{ is optimal}\} \supset E_0 \cap E_1 \cap E_2$, where $E_0 = \{\tilde{M} \text{ is true}\}$.

7.3. *Proof of Theorem* 2. First note that like the fence procedure, the F-B fence is guaranteed to stop at some point. This is because, otherwise, one keeps adding the parameters until one gets the full model, which automatically satisfies the fence inequality (note that in this case $\tilde{M}$ is chosen as the full model).

Next we show that w.p. $\to 1$, $M_0^\dagger$ is a true model. Suppose that this is not the case. Then there is an incorrect model, say, $M$, such that

(10) $$\mathrm{P}(M_0^\dagger = M) \geq \delta,$$

where $\delta > 0$ is a constant. Since $\tilde{M}$ is a true model, we have by the proof of Theorem 1 that w.p. $\to 1$, $\hat{Q}_M > \hat{Q}_{\tilde{M}} + c_n \hat{\sigma}_{M,\tilde{M}}$. On the other hand, $M_0^\dagger = M$ implies that $\hat{Q}_M \leq \hat{Q}_{\tilde{M}} + c_n \hat{\sigma}_{M,\tilde{M}}$ [because $M_0^\dagger$ has to satisfy (3)]. Thus, we have $\mathrm{P}(M_0^\dagger = M) \leq \mathrm{P}(\hat{Q}_M \leq \hat{Q}_{\tilde{M}} + c_n \hat{\sigma}_{M,\tilde{M}}) \to 0$, which contradicts (10).

We next show that w.p. $\to 1$, no proper submodel of $M_0^\dagger$ is a true model. Suppose that this is not true. Then there is a true model $M_1$ and a constant $\delta > 0$ such that $\mathrm{P}(M_1 \subset M_0^\dagger) \geq \delta$. Hereafter, the notation $M_1 \subseteq M_2$ ($M_1 \subset M_2$) means that $M_1$ is a (proper) submodel of $M_2$. Suppose that under $M_0^\dagger$, $X\beta + Z\alpha = \sum_{r \in R_0} X_r \beta_r + \sum_{s \in S_0} Z_s \alpha_s$, and, under $M_1$, the same expression holds with $R_0, S_0$ replaced by $R_1, S_1$, respectively. Define $R_{10} = R_1 \cup \{r_1, \ldots, r_{a-1}\}$, $S_{10} = S_0$, if $R_1 \subset R_0$, $S_1 \subseteq S_0$ and $R_0 \setminus R_1 = \{r_1, \ldots, r_a\}$; $R_{10} = R_0$, $S_{10} = S_1 \cup \{s_1, \ldots, s_{b-1}\}$, if $R_1 = R_0$, $S_1 \subset S_0$ and $S_0 \setminus S_1 = \{s_1, \ldots, s_b\}$; and $R_{10} = R_1$, $S_{10} = S_1$ otherwise. Let $M_{10}$ be the model corresponding to $R_{10}$ and $S_{10}$. Then $M_1 \subset M_0^\dagger$ implies that $M_{10} \subset M_0^\dagger$ with one less parameter, hence we must have $\hat{Q}_{M_{10}} > \hat{Q}_{\tilde{M}} + c_n \hat{\sigma}_{M_{10},\tilde{M}}$ by the definition of $M_0^\dagger$. It follows that

(11) $$\mathrm{P}(\hat{Q}_{M_{10}} > \hat{Q}_{\tilde{M}} + c_n \hat{\sigma}_{M_{10},\tilde{M}}) \geq \delta.$$

On the other hand, we have by the proof of Theorem 1 that for any true model $M$, w.p. $\to 1$, $\hat{Q}_M \leq \hat{Q}_{\tilde{M}} + c_n \hat{\sigma}_{M,\tilde{M}}$. Since $M_{10}$ is always a true model, it follows that $\mathrm{P}(\hat{Q}_{M_{10}} > \hat{Q}_{\tilde{M}} + c_n \hat{\sigma}_{M_{10},\tilde{M}}) \leq \sum_{M \text{ true}} \mathrm{P}(\hat{Q}_M > \hat{Q}_{\tilde{M}} + c_n \hat{\sigma}_{M,\tilde{M}}) \to 0$, which contradicts (11).

7.4. *Proof of Theorem* 3. (i) For any $\varepsilon$, $\eta > 0$, let $\delta$, $x_{n,j}$, $j = 1, 2, 3$, $N$ and $N^*$ be as in Assumptions A7 and A8. Then when $n \geq N$ and $n^* \geq N^*$, the following arguments hold with probability $> 1 - \eta$.

For $j = 1, 3$, we have $P^*(x_{n,j}) > P(x_{n,j}) - \delta > 1 - 5\delta \geq 1/2$. It follows that $p^*(x_{n,j}) = \max_{M \in \mathcal{M}} P^*(M_0(x_{n,j}) = M) = P^*(x_{n,j}) < P(x_{n,j}) + \delta \leq 1 - 2\delta$. Similarly, $p^*(x_{n,2}) = P^*(x_{n,2}) > P(x_{n,2}) - \delta > 1 - 2\delta$. Thus, there is $c_n^* \in (x_{n,1}, x_{n,3})$ which is the maximum of $p^*$ over $[x_{n,1}, x_{n,3}]$. Furthermore, we have $p^*(c_n^*) \geq p^*(x_{n,2}) > 1 - 2\delta$.

(ii) If $M_{\text{opt}} = M_{\text{f}}$, then $q_n \sim a_n/b_n$, hence $q_n^{-1} = (b_n/a_n)O(1) = o(1)$. Also, by Assumption A8, for any $\eta > 0$, there is $L > 0$ such that $P(q_n/q_n^* > L) < \eta$. Choose $N_1 \geq 1$ such that $q_n^{-1} < 1/L$ when $n \geq N_1$. Then, when $n \geq N_1$, we have, w. p. $> 1 - \eta$, $(q_n^*)^{-1} = q_n^{-1}(q_n/q_n^*) < 1$, hence $q_n^* > 1$, hence $c_n^* = 0$. On the other hand, by Assumption A7, there is $N_2 \geq 1$ such that $P(\min_{M \in \mathcal{M}, M \neq M_{\text{f}}} \hat{Q}_M > \hat{Q}_{M_{\text{f}}}) > 1 - \eta$, if $n \geq N_2$. Let $N = N_1 \vee N_2$, then $P(M_0^* = M_{\text{f}}) > 1 - 2\eta$, if $n \geq N$.

If $M_{\text{opt}} = M_*$, then by similar arguments, it can be shown that $r_n^* = o_P(1)$ and $q_n^* = o_P(1)$. Thus, for any $\eta > 0$, there is $N \geq 1$ such that when $n \geq N$ we have, w.p. $> 1 - \eta$, $q_n^* \leq 1$ and $r_n^* < 1$, hence $c_n^* = B^*$, hence $M_0^* = M_*$.

If $M_{\text{opt}} \notin \{M_{\text{f}}, M_*\}$, note that $\{M_0^* = M_{\text{opt}}\} \supset \{r_{\text{opt}} \leq c_n^*, r_{\text{w} \leq} > c_n^*\} \supset \{r_{\text{opt}} \leq x_{n,1}, r_{\text{w} \leq} > x_{n,3}\}$, if $c_n^* \in (x_{n,1}, x_{n,3})$. Therefore, by (i), for any $\varepsilon$, $\eta > 0$, we have

$$\begin{aligned}
P(M_0^* = M_{\text{opt}}) &\geq P(M_0^* = M_{\text{opt}}, c_n^* \in (x_{n,1}, x_{n,3})) \\
&\geq P(r_{\text{opt}} \leq x_{n,1}, r_{\text{w} \leq} > x_{n,3}, c_n^* \in (x_{n,1}, x_{n,3})) \\
&\geq F_{\text{opt}}(x_{n,1}) - F_{\text{w} \leq}(x_{n,3}) - P(c_n^* \notin (x_{n,1}, x_{n,3})) \\
&> 1 - 2\varepsilon - \eta, n \geq N, n^* \geq N^*.
\end{aligned}$$

**Acknowledgments.** The authors are grateful to the Associate Editor and referees for their constructive comments that have led to the development of adaptive fence and other improvements.

J. Jiang  
T. Nguyen  
Department of Statistics  
University of California, Davis  
One Shields Ave  
Davis, California 95616  
USA  
E-mail: jiang@wald.ucdavis.edu  
    tnguyen@wald.ucdavis.edu  

Z. Gu  
ALZA Corporation  
1900 Charleston Road, M11-4  
Mountain View, California 94039  
USA  
E-mail: zgu26@alzus.jnj.com  

J. Sunil Rao  
Department of Epidemiology  
  and Biostatistics  
Case Western Reserve University  
10900 Euclid Ave  
Cleveland, Ohio 44106  
USA  
E-mail: sunil@hal.cwru.edu